\documentclass[12pt]{article}

\usepackage[T2A]{fontenc}
\usepackage[utf8]{inputenc}
\usepackage[russian]{babel}
\usepackage{amssymb,amsfonts,amsthm,amsmath}
 \usepackage[font=onehalfspacing,margin=1cm]{caption}
\usepackage[a4paper,lmargin=2cm,margin=2cm,rmargin=2cm]{geometry}

\usepackage[doublespacing]{setspace}
\usepackage{wrapfig}

\newtheorem{prep}{Утверждение}
\renewcommand{\Pr}{{\mathbf P}}
\usepackage{tikz}
\usetikzlibrary{calc}

\title{Меры Хаусдорфа в задаче Бертрана \\о случайной хорде}
\author{А. В. Зорин\footnote{Национальный исследовательский Нижегородский
    государственный университет им. Н. И. Лобачевского, г. Нижний Новгород,
    Россия; e-mail: andrei.zorine@itmm.unn.ru}} \date{}
\begin{document}

\maketitle

\begin{abstract}
  Множество хорд некоторой окружности заданного радиуса представляется
  как метрическое пространство относительно метрики, введенной
  Хаусдорфом. Устанавливается вид открытых и замкнутых шаров относительно этой
  метрики. Рассматривается семейство внешних мер Хаусдорфа, порожденных это
  метрикой.  Вычисляется размерность Хаусдорфа открытых и замкнутых
  шаров. Вводится аналог непрерывного равномерного распределения и предлагается
  новое решение задачи Бертрана со старым ответом.
  
  \emph{Ключевые слова:}\/ метрическое пространство, внешняя мера, мера
  Хаусдорфа, парадокс Бертрана
\end{abstract}

\bigskip

\textbf{1. История парадокса Бертрана и построение вероятностных пространств.}

Следующая задача Ж.~Бертрана хорошо известна: \textit{<<Найти вероятность, с
  которой случайно выбранная хорда на окружности окажется длинее стороны
  равностороннего треугольника, вписанного в эту окружность>>}. Сам
Бертран~\cite{Bertrand} приводит три решения этой задачи, дающие разные ответы
(от одной черверти до одной второй). Все решения заключаются в том, что <<выбор
хорды наудачу>> сводится к выбору наудачу точки --- на некотором отрезке, в
некотором круге --- и в предположении о равновозможности всех положений такой
точки в выборочном пространстве. По мнению Бертрана, затруднение вызвано
попыткой применить понятие выбора наудачу для элемента \emph{бесконечноного}
множества.  В настоящее время общепризнано
\cite{Gnedenko,Sekely,Fedotkin,Gantsevich,Kaushik:2022}, что проблема
заключается в недоопределенном механизме выбора хорды и что разные решения
предполагают равномерное распределние у разных \emph{параметризаций}
задачи. Предлагаются новые решения, с новым ответом. Делается попытка выделить
<<верное>> решение путем наложения дополнительного требования инвариантности
распределения параметров относительно изменения масштаба и сдвигов. Чисто
умозрительные требования типа <<инвариантности>>, кажется, еще дальше уводят от
вопроса: как физические условия проведения эксперимента определяют вероятность
события? В книге~\cite[с.106--109]{Fedotkin} задача о хорде интерпретируется как
задача о трещине на гончарном круге, делается попытка сопоставить три разные
решения (ответа) с разными способами производства и материалом изделия. При
таком подходе все три решения оказываются физически реализуемыми, и вопрос о
том, какое из них <<верное>>, сам лишается смысла. Правдоподобность трех
случайных механизмов подтверждается свидетельствами археологии о ценообразовании
у гончарных кругов в Древней Греции (по сообщению цитируемой
книги~\cite{Fedotkin}).

Особо отметим здесь работу~\cite{Gantsevich}, в которой задача решается путем
предельного перехода от конечной задачи к бесконечной: хорды представлены как
диагонали вписанного \emph{правильного} многоугольника. Затем число вершин
многоугольника устремляется к бесконечности. Данная схема предельного перехода
может дать только бесконечное счетное подмножество из континуального множества
всех хорд данной окружности, причем хорд, расположенных довольно
<<регулярно>>. Тем не менее, данная статья интересна тем, что вместо работы с
параметрическим представлением хорды объектом случайного выбора снова становится
сама хорда, как изначально задумывал Бертран. В рамках такого
непараметрического, геометрического подхода, в работе \cite{Aldou} решается
задача о <<случайной триангуляции окружности>>. Для случайной триангуляции также
можно исследовать распределение длины, например, самой длинной хорды. Метод
задания распределения вероятностей на множестве всех триангуляций в цитируемой
работе отличается от обсуждавшихся выше и отличается от предлагаемого в
настоящей заметке.

Исключительное положение непрерывного равномерного распределения среди всех
непрерывных распределений вызвано тем, что оно порождено простым геометрическим
понятием длины, площади, объема.  Следовательно, чтобы решить задачу о хорде
геометрически, надо ввести понятие геометрической меры для любого множества
хорд. Это возможно сделать, если превратить множество хорд в метрическое
пространство и воспользоваться конструкцией \emph{внешней меры}
Хаусдорфа~\cite{Edgar} на абстрактном метрическом пространстве. Именно
реализации этой идеи посвящена данная работа.\sloppy

\bigskip

\textbf{2. Множество хорд как метрическое пространство.}
Пусть $X$~--- множество всех хорд заданной окружности
радиуса~$R>0$. Договоримся, что точки окружности не являются хордами
(вырожденными хордами с совпадающими концами).  Множество~$X$ имеет мощность
континуума. Его можно превратить в метрическое пространство, введя расстояние
между произвольными хордами $\chi\in X$, $\chi'\in X$ как расстояние Хаусдорфа~\cite{Hausdorf} между двумя
плоскими компактными множествами:
\[
  \rho(\chi, \chi') = \max\raisebox{-0.25ex}{$\mathsurround=0pt \bigl\{$}
  \max_{P\in \chi} \min_{Q\in\chi'}|PQ|,\;
  \max_{P\in \chi'} \min_{Q\in\chi}|PQ|
  \raisebox{-0.25ex}{$\mathsurround=0pt \bigr\}$}.
\]
Читатель может легко проверить, что расстояния Хасдорфа между хордами на 
рис.~\ref{AZ:fig:1} равны длине отрезка $BD$. Например, слева мы видим две
<<сходящиеся>>, но не пересекающиеся хорды. Дуга $DB$ больше дуги $AC$, самая
удаленная точка хорды $AB$ от хорды $CD$ есть точка $B$, расстояние --- длина
отрезка $BD$, а самая удаленная точка хорды $CD$ от хорды $AB$ есть точка $D$,
расстояние --- длина перпендикуляра, опущенного из точки $D$ на хорду $AB$. По
свойству перпундикуляра, его длина меньше наклонного отрезка $BD$, значит,
расстояние Хаусдорфа между этими двумя хордами есть длина отрезка $BD$.

\begin{figure}[h]
  \qquad
  \begin{tikzpicture}
  \def\ptsize{2.0pt};
  \draw[thick] (0,0) circle (1.7cm);
  \coordinate[label=above right:\(A\)] (A) at (30:1.7cm); 
  \coordinate[label=below left:\(B\)] (B) at (180:1.7cm);
  \coordinate[label=above right:\(C\)] (C) at (5:1.7cm);
  \coordinate[label=below left:\(D\)] (D) at (230:1.7cm);
  \draw (A) -- (B);
  \draw (C) -- (D);
  \draw[black!50,thick] (D) -- ($(B)!(D)!(A)$);
  \draw[black!50,thick] (B)--(D);
  \fill[fill=white,draw=black,thick] (A) circle (2pt);
  \fill[fill=white,draw=black,thick] (B) circle (2pt);
  \fill[fill=white,draw=black,thick] (C) circle (2pt);
  \fill[fill=white,draw=black,thick] (D) circle (2pt);
\end{tikzpicture}
\quad
\begin{tikzpicture}
  \def\ptsize{2.0pt};
  \draw[thick] (0,0) circle (1.7cm);
  \coordinate[label=above right:\(A\)] (A) at (30:1.7cm);
  \coordinate[label=below left:\(B\)] (B) at (180:1.7cm);
  \coordinate[label=above right:\(C\)] (C) at (5:1.7cm);
  \coordinate[label=below left:\(D\)] (D) at (230:1.7cm);
  \draw (A) -- (D);
  \draw (C) -- (B);
  \draw[black!50] (B) -- ($(D)!(B)!(A)$);
  \draw[black!50] (B)--(D);
  \fill[fill=white,draw=black,thick] (A) circle (2pt);
  \fill[fill=white,draw=black,thick] (B) circle (2pt);
  \fill[fill=white,draw=black,thick] (C) circle (2pt);
  \fill[fill=white,draw=black,thick] (D) circle (2pt);
\end{tikzpicture}
\quad
\begin{tikzpicture}
  \def\ptsize{2.0pt};
  \draw[thick] (0,0) circle (1.7cm);
  \coordinate[label=above right:\(A\)] (A) at (30:1.7cm);
  \coordinate[label=below left:\(B\)] (B) at (150:1.7cm);
  \coordinate[label=above right:\(C\)] (C) at (5:1.7cm);
  \coordinate[label=below left:\(D\)] (D) at (175:1.7cm);
  \draw (A) -- (B);
  \draw (C) -- (D);
  \fill[fill=white,draw=black,thick] (A) circle (2pt);
  \fill[fill=white,draw=black,thick] (B) circle (2pt);
  \fill[fill=white,draw=black,thick] (C) circle (2pt);
  \fill[fill=white,draw=black,thick] (D) circle (2pt);
\end{tikzpicture}
\caption{\label{AZ:fig:1}%
  Примеры на вычисление расстояния между хордами
}
\end{figure}

Пусть $\varepsilon$~--- малое положительное число. Выясним, как выглядит
замкнутый шар
$$
\bar{\mathcal B}(\chi, \varepsilon)=\{\chi'\in X\colon \rho(\chi,
\chi')\leqslant \varepsilon\}
$$ радиуса $\varepsilon$ с
центром в произвольной хорде $\chi\in X$. Хорда $\chi'$ входит в
$\varepsilon$-окрестность хорды $\chi$, если концы хорды $\chi'$ лежат от концов
хорды $\chi$ на расстоянии не больше $\varepsilon$. На рис.~\ref{AZ:fig:ball}
это всякая хорда, один конец которой лежит на дуге $AC$, а другой --- на дуге
$BD$. Открытый шар
$$
{\mathcal B}(\chi, \varepsilon)=\{\chi'\in X\colon \rho(\chi,
\chi')< \varepsilon\}
$$ радиуса $\varepsilon$ и с
центром в произвольной хорде $\chi\in X$ состоит из хорд $\chi'$, концы которых
не совпадают ни с одной из точек $A$, $B$, $C$, $D$. Граница шара $\bar{\mathcal
  B}(\chi, \varepsilon)$  состоит из  хорд, один конец которых совпадает с одной
из точек $A$, $B$, $C$, $D$, а другой располагается произвольно на
противоположной дуге (например, точка $A$ и произвольная точка дуги $BD$).

\begin{figure}[h]
  \centering
  \begin{tikzpicture}
  \def\ptsize{2.0pt};
  \draw[thick] (3.8637, 1.0352) arc (15:165:4cm);
  \coordinate[label=above right:\(A\)] (A) at (50:4cm);
  \coordinate[label=above left:\(B\)] (B) at (130:4cm);
  \coordinate[label=above right:\(C\)] (C) at (20:4cm);
  \coordinate[label=above left:\(D\)] (D) at (160:4cm);
  \coordinate[label=above right:\(C\)] (C) at (20:4cm);
  \coordinate[label=above left:\(D\)] (D) at (160:4cm);
  \coordinate(A1) at (35:4cm);
  \coordinate(B1) at (145:4cm);
  \draw[thin] (A) -- (B);
  \draw[thick] (A1) -- (B1);
  \node  at (27.5:3.8cm) {$\varepsilon$};
  \node  at (42.5:3.8cm) {$\varepsilon$};
  \node  at (152.5:3.8cm) {$\varepsilon$};
  \node  at (137.5:3.8cm) {$\varepsilon$};
  \node  at (0,2.5cm) {$\chi$};
  \draw[thin] (C) -- (D);
  \draw[black!50] (B)--(B1)--(D);
  \draw[black!50] (A)--(A1)--(C);
  \fill[fill=white,draw=black,thick] (A1) circle (2pt);
  \fill[fill=white,draw=black,thick] (B1) circle (2pt);
  \fill[fill=white,draw=black,thick] (A) circle (2pt);
  \fill[fill=white,draw=black,thick] (B) circle (2pt);
  \fill[fill=white,draw=black,thick] (C) circle (2pt);
  \fill[fill=white,draw=black,thick] (D) circle (2pt);
\end{tikzpicture}
\caption{\label{AZ:fig:ball}%
К построению $\varepsilon$-окрестности хорды $\chi$}
\end{figure}

Обратно, для двух параллельных хорд $AB$, $CD$, множество всех хорд, один конец
которых лежит на дуге $AC$, а другой на дуге $BD$, является замкнутым шаром с
центром в хорде, концы которой лежат на серединах дуг $AC$ и $BD$.  Нам будет
удобно  называть такой замкнутый шар \emph{замкнутой трубкой} между
параллельными хордами $AB$ и $CD$, а открытый шар с тем же центром и радиусом
--- \emph{открытой трубкой}. Обозначим открутую и замкнутую трубки между хордами
$AB$ и $CD$ через $\mathcal T(AB, CD)$ и $\bar{\mathcal T}(AB, CD)$
соответственно. В дальнейшем удобно предполагать, что по окружности точки вседа
лежат в порядке $A$, $B$, $D$, $C$ при обходе в определенном направлении, так
что дуги между хордами всегда $AC$ и $BD$.

Диаметр множества $A\subset X$ определяется равенством
$d(A)=\sup_{\chi, \chi' \in A} \rho(\chi, \chi').$

\begin{prep}
  Диаметр $d(X)$ пространства $X$ равен $2R$.
\end{prep}

\begin{proof}
  Очевидно, $\rho(\chi, \chi')\leqslant 2R$ из геометрических соображений (любой
  отрезок $PQ$ с концами на разных хордах, который дает искомое расстояние
  Хаусдорфа, целиком лежит в круге и может быть дополнен до некоторой хорды,
  которая не превосходит диаметра). С другой стороны, рассмотрим
  последовательность $(\chi_n, \chi'_n)$, $n=1$, $2$, \ldots{} пар хорд, так что
  в каждой паре хорды параллельны и расстояние от точек пересечения хорд с общим
  перпендикулярным им диаметром до центра круга равно $(R-1/n)$. Очевидно,
  \[
    d(X)\geqslant \rho(\chi_n, \chi_n')=2R-2/n\to2.
  \]
  Доказательство завершено.
\end{proof}

\begin{prep}
  Пусть параллельные хорды $\chi$, $\chi'$ заключают между собой дуги длины
  $\gamma$, расположенные на одной полуокружности (т.е. обе хорды расположены по
  одну сторону от параллельного им диаметра). Тогда диаметры трубок
  $\mathcal T(\chi, \chi')$ и $\bar{\mathcal T}(\chi, \chi')$ равны
  $2R\sin\dfrac{\gamma}{2R}$.
\end{prep}

\begin{figure}[h]
   \centering
 \begin{tikzpicture}
   \def\ptsize{2.0pt};
   \coordinate[label=above:$O$] (O) at (0,0);
   \coordinate[label=above right:\(A\)] (A) at (45:2cm);
   \coordinate[label=above right:\(B\)] (B) at (-15:2cm);
   \coordinate[label=above right:\(C\)] (C) at (90:2cm);
   \coordinate[label=below right:\(D\)] (D) at (-60:2cm);
   \coordinate[label=below right:\(F\)] (F) at (-37.5:2cm);
   \draw[thick] (0,0) circle (2cm);
   \draw[thick] (A)--(B);
   \draw[thick] (C)--(D);
   \draw[black!50] (O)--(A);
   \draw[black!50] (O)--(C);
   \draw[black!50] (A)--(C);
   \fill (O) circle(2pt);
   \foreach \point in {A,B,C,D}
   \fill[fill=white,draw=black,thick] (\point) circle (2pt);
   \node at (67:2.2cm) {$\gamma$};
 \end{tikzpicture}
 \caption{\label{AZ:fig:tube-diam} К вычислению диаметра трубки}
\end{figure}

\begin{proof}
  Даны две параллельные хорды $\chi=AB$ и $\chi'=CD$. Угол $COA$ в радианах
  равен $\gamma/R$, тогда длина отрезка $AC$ равна $2R\sin(\gamma/2R)$
  (см. рис. \ref{AZ:fig:tube-diam}).
\end{proof}

\begin{prep}
  \label{AZ:prep:compact}
  Из всякого счетного покрытия замкнутой трубки $\bar{\mathcal T}(AB, CD)$
  открытыми трубками ${\mathcal T}(A_nB_n, C_nD_n)$, $n=1$, $2$, \ldots{} можно
  выбрать конечное покрытие.
\end{prep}

\begin{proof}
  Доказательство проводится через деление дуг $AC$ и $BD$ пополам, как при
  доказательстве теоремы Гейне--Бореля для прямоугольника~\cite{Zorich}.
\end{proof}

\bigskip

\textbf{3. Мера Хаусдорфа и размерность Хаусдорфа для некоторых множеств хорд.}
На метрическом пространстве $(X, \rho)$ для произвольного положительного
действительного числа $s$ можно определить \emph{$s$-мерную внешнюю меру
  Хаусдорфа} $\bar{\mathcal H}^s(\cdot)$ следующим
образом~\cite{Edgar}. Определяется для всякого $\varepsilon>0$ внешняя мера
\emph{по методу I}
$$
\bar{\mathcal H}^s_\varepsilon (E)= \inf\Bigl\{ \sum_{i=1}^\infty (d(E_i))^s\colon
E\subset \bigcup_{i=1}^\infty E_i, \; d(E_i)\leqslant \varepsilon,\; i=1, 2,
\ldots \Bigr\},
$$
где $\inf$ берется по всем счетным покрытиям множества $E$ открытыми или
закрытыми трубками $E_i$ диаметра не больше $\varepsilon$. Далее, внешняя мера
Хаусдорфа получается так называемем \emph{методом II} путем предельного перехода
$$
\bar{\mathcal H}^s(E) = \lim_{\varepsilon\to0} \bar{\mathcal
  H}^s_\varepsilon(E)=\sup_{\varepsilon>0} \bar{\mathcal
  H}^s_\varepsilon(E) .
$$
Данная конструкция гарантирует~\cite{Edgar} измеримость всех открытых подмножеств множества
$X$. Из теории меры Хаусдорфа известно, что существует такое число $s_0$, что
$\bar{\mathcal H}^{s_0}(E)=\infty$ для всех $s<s_0$, $\bar{\mathcal H}^{s_0}(E)=0$ для всех
$s>s_0$. Это число $s_0$ называется хаусдорфовой размерностью множества $E$.

\begin{prep}
  \label{AZ:prep:tube}
  Для произвольной замкнутой трубки $\bar{\mathcal T}(AB,CD)$, не содержащей
  диаметра окружности и опирающейся на дуги длины $\gamma$, ее внешняя мера
  Хаусдорфа $\bar{\mathcal H}^s(\bar{\mathcal T})$ равна $0$ при $s>2$, равна
  $\gamma^2$ при $s=2$ и равна $\infty$ при $s<2$. Хаусдорфова размерность
  произвольной трубки с $\gamma>0$ равна двум.
\end{prep}
\medskip

\begin{proof}
  Поскольку любая замкнутая трубка компактна (Утв.~\ref{AZ:prep:compact}),
  достаточно рассматривать только конечные покрытия. Пусть $n$ определяется
  условием: $(n-1)\varepsilon<\gamma\leqslant n\varepsilon$,
  т.е. $n=\gamma/\varepsilon$, если это отношение является целым числом, в
  противном случае $n=[\gamma/\varepsilon]+1$. Разобьем дугу $AC$ точками $E_1$,
  $E_2$, \ldots, $E_{n-1}$, а дугу $BD$ точками $F_1$, $F_2$, \ldots, $F_{n-1}$
  на $n$ равных дуг длины $\gamma/n\leqslant\varepsilon$. Обозначим $E_0=A$,
  $E_n=C$, $F_0=B$, $F_n=D$. Тогда набор трубок
  $\bar{\mathcal T}(E_{k-1}E_k, F_{l-1},F_l)$, $k,l=\overline{1,n}$ образует
  конечное покрытие и выполнено неравенство
  \begin{equation}
    \bar{\mathcal H}^s_\varepsilon (\bar{\mathcal T}(AB,CD) )\leqslant
    \sum_{k=1}^n \sum_{l=1}^n (d(\mathcal T(E_{k-1}E_k, F_{l-1},F_l)))^s=
    n^2 \Bigl(2R\sin\bigl( \frac{\gamma}{2Rn}\bigr)\Bigr)^s.
    \label{AZ:eq:Hsup}
  \end{equation}
При $s>2$ правая часть стремится к нулю при  $\varepsilon\to\infty$ до
  нуля. Первая часть утверждения доказана.

  Для доказательства второй и третьей части утверждения положим $s=2$ и оценим
  $\bar{\mathcal H}^2_\varepsilon(\bar{\mathcal T}(AB,CD))$ снизу. Каждую хорду
  $EF$ из трубки $\bar{\mathcal T}(AB,CD)$ с вершиной $E$ на $AC$ и вершиной $F$
  на $BD$ можно изобразить как точку квадрата со стороной $\gamma$, откладывая
  по оси абсцисс длину дуги $AE$, а по оси ординат длину дуги $BF$
  (рис.~\ref{AZ:fig:sq}). При этом, каждая трубка из покрытия будет изображаться
  квадратом, а все конечное покрытие трубки будет изображаться покрытием
  квадрата со стороной $\gamma$ конечным числом $N$ квадратов со сторонами
  $x_1$, $x_2$, \ldots, $x_N$ длины не больше $\varepsilon$.
  
  \begin{figure}
    \centering
    \begin{tikzpicture}[scale=0.67]
      \def\ptsize{2.0pt}; \coordinate[label=above:$O$] (O) at (0,0);
      \coordinate[label=above right:\(A\)] (A) at (45:2cm);
      \coordinate[label=above right:\(B\)] (B) at (-15:2cm);
      \coordinate[label=above left:\(C\)] (C) at (90:2cm);
      \coordinate[label=below left:\(D\)] (D) at (-60:2cm);
      \draw[thick] (0,0) circle (2cm);
      \draw[thick] (A)--(B);
      \draw[thick] (C)--(D);
      \draw[black!50,thick]  (50:2cm)--(-50:2cm);
      \draw[black!50,thick]  (55:2cm)--(-55:2cm);
      \fill (O) circle(2pt); \foreach \point in {A,B,C,D}
      \fill[fill=white,draw=black,thick] (\point) circle (2pt); \node at
      (67:2.2cm) {$\gamma$};
    \end{tikzpicture}
    \qquad \raisebox{1.65cm}{$\longrightarrow$} \qquad
    \begin{tikzpicture}
      \draw[->,thick] (-0.2,0) -- (3.2,0);
      \draw[->,thick] (0,-0.2) -- (0,3.2);
      \draw (2.7,0)--(2.7,2.7)--(0,2.7);
      \fill[black!50] (0.3,2.1)--(0.6,2.1)--(0.6,2.4)--(0.3,2.4)--cycle;
      \draw (-0.1,-0.1)--(0.5,-0.1)--(0.5,0.5)--(-0.1,0.5)--cycle;
      \draw (0.5,-0.1)--(0.5,1.5)--(2.1,1.5)--(2.1,-0.1)--cycle;
      \draw (0.4,0.4)--(0.4,2.8)--(2.8,2.8)--(2.8,0.4)--cycle;
      \node [anchor=north] at (2.7,0) {$\gamma$};
      \node [anchor=east] at (0,2.7) {$\gamma$};
      \node [anchor=east] at (-0.1,0.2) {$x_1$};
      \node [anchor=east] at (2.1,1.2) {$x_2$};
      \node [anchor=south] at (1.5,2.7) {$x_3$};
      \node [anchor=west] at (0.45,2.4) {$x_4$};
    \end{tikzpicture}
    \caption{\label{AZ:fig:sq}%
      Покрытие трубками и отображение в квадраты. Трубке между серыми хордами
      слева отвечает серый квадрат справа. У сторон квадратов указаны их размеры
      ($x_1$, $x_2$, \ldots{})}
  \end{figure}
  
  Суммарная площадь этих квадратов равна
  \begin{equation}
    x_1^2+x_2^2+\ldots+x_N^2\geqslant \gamma^2.\label{AZ:eq:1}
  \end{equation}
  С другой стороны,
  \begin{equation}
    0\leqslant x_i \leqslant \varepsilon, \qquad i=1,2,\ldots,N,
    \label{AZ:eq:2}
  \end{equation}
  откуда
  $$    x_1^2+x_2^2+\ldots+x_N^2\leqslant N\cdot \varepsilon^2.$$
  Значит, число $N$ покрывающих трубок должно удовлетворять неравенству
  \begin{equation}
    N\geqslant(\gamma^2/\varepsilon^2).
    \label{AZ:eq:3}
  \end{equation}
  Требуется найти точную нижнюю грань величины
  $$
  f_N(x_1, x_2, \ldots, x_N)=\sum_{i=1}^N \Bigl(
     2R \sin\Bigl(\frac{x_i}{2R}\Bigr)
  \Bigr)^2
  $$
  по всем допустимым $N$, $x_1$, \ldots, $x_N$, удовлетворяющим неравенствам
  \eqref{AZ:eq:1}, \eqref{AZ:eq:2} и \eqref{AZ:eq:3}. Решение задачи
  оптимизации $$f_N(x_1, \ldots, x_N)\to\min_{N,x_1,\ldots,x_N}$$ при указанных
  условиях не обязано соответствовать какому-то покрытию, но дает нижнюю оценку
  величины $\bar{\mathcal H}^s(\bar{\mathcal T}(AB,CD))$.
  
  Заменим каждый из синусов многочленом Тейлора с остатком в форме Лагранжа,
  получим:
  $$
  \sum_{i=1}^N \Bigl(2R \sin\Bigl(\frac{x_i}{2R}\Bigr)\Bigr)^2=
  \sum_{i=1}^N\Bigl(x_i-\cos\Bigl(\frac{ \theta_i x_i}{2R} \Bigr)
  \frac{x_i^3}{24R^2}
  \Bigr)^2,
  $$
  где $0\leqslant  \theta_i \leqslant 1$ для всех $i=1$, $2$, \ldots, $N$.
  Далее, имеют место оценки
  \begin{multline}
    \sum_{i=1}^N\Bigl(x_i-\cos\Bigl(\frac{ \theta_i x_i}{2R} \Bigr)
    \frac{x_i^3}{24R^2} \Bigr)^2=  \sum_{i=1}^N x_i^2 -\sum_{i=1}^N
    \cos\Bigl(\frac{ \theta_i x_i}{2R} \Bigr) \frac{x_i^4}{12R^2}+
    \sum_{i=1}^N \cos^2\Bigl(\frac{ \theta_i x_i}{2R} \Bigr)
    \frac{x_i^6}{576R^4}\geqslant
    \\ \geqslant \sum_{i=1}^N x_i^2
    -\sum_{i=1}^N \frac{x_i^4}{12R^2} \geqslant \sum_{i=1}^N x_i^2
    -\sum_{i=1}^N \frac{\varepsilon^2x_i^2}{12R^2} =
    \Bigl(1-\frac{\varepsilon^2}{12R^2}\Bigr)\sum_{i=1}^N x_i^2 \geqslant
    \Bigl(1-\frac{\varepsilon^2}{12R^2}\Bigr) \gamma^2.
    \label{AZ:eq:H2low}
  \end{multline}
  Теперь из \eqref{AZ:eq:Hsup} и
  \eqref{AZ:eq:H2low} заключаем, что
  $$
  \lim_{\varepsilon\to0} \bar{\mathcal H}^2_\varepsilon( \bar{\mathcal
    H}(AB,CD))=\gamma^2.
  $$
  Наконец, из свойств внешней меры Хаусдорфа следует, что
  $\bar{\mathcal H}^s( \bar{\mathcal T}(AB,CD))=\infty$ при $s<2$
  и что размерность Хаусдорфа трубки $\bar{\mathcal T}(AB,CD)$ равна двум.
\end{proof}

Таким образом, мы вычислили внешнюю меру трубки на основании определения
внешней меры Хаусдорфа. Следует отметить, что возникшая промежуточная
конструкция покрытия квадрата квадратами не полностью аналогична покрытию
прямоугольников прямоугольниками из теории плоской меры
Лебега~\cite{Kolmogorov:Fomin}. Дело в том, что конечные объединения
квадратов (со сторонами, параллельными координатным осям) не образуют
теоретико-множественного \emph{полукольца}. Можно так расположить квадраты с
нахлестом, что потребуется не менее чем бесконечное счетное покрытие
квадратами, тогда как хватило бы конечного числа прямоугольников, правда, с
несоизмеримыми сторонами.

В оставшейся части работы сконцентрируемся на случае хаусдорфовой размерности,
то есть $s=2$. Теперь мы можем вычислять меры произвольных измеримых множеств,
используя общие теоремы и свойства мер, такие как свойство счетной аддитивности,
свойство непрерывности для монотонных классов множеств.

\begin{prep}
  Мера Хаусдорфа множества хорд, заключенных между двумя непересекающимся
  хордами, дуги между которыми равны $\gamma_1$, $\gamma_2$, равна
  $\gamma_1\gamma_2$.  
\end{prep}
\begin{proof}
  Выберем малое $\varepsilon>0$. Возьмем $k_i$ из условия
  $(k_i-1)\varepsilon<\gamma_i\leqslant k_i \varepsilon$, $i=1$, $2$. Тогда дугу
  $\gamma_i$ можно покрыть $k_i$ дугами длины~$\varepsilon$. Значит, всё
  рассматриваемое множество хорд можно покрыть $k_1\cdot k_2$ трубками меры
  $\varepsilon^2$. Переходя к пределу, имеем
  $$
  \lim_{\varepsilon\to0}
   (k_1\cdot k_2)
  \varepsilon^2=\gamma_1\gamma_2. \qedhere
  $$
\end{proof}

\begin{prep}
  \label{AZ:prep:1}
  Пусть хорда $AB$ стягивает дугу длины $\gamma$. Множество хорд, концы которых
  лежат на этой дуге, имеет меру $\gamma^2/2$.
\end{prep}

\begin{proof}
  Пусть $C_1$~--- середина дуги, стягиваемой хордой $AB$. Все хорды из
  рассматриваемого множества принадлежат одному из трех непересекающихся
  множеств: 1)~один конец на дуге $AC_1$, другой конец на дуге $C_1B$
  (<<вырожденная трубка>>, одна из ограничивающих хорд вырождается в точку); 2) оба
  конца лежат на дуге $AC_1$ 3) оба конца лежат на дуге $C_1B$.

  \begin{figure}[h]
    \centering
    \begin{tikzpicture}[scale=0.75]
      \def\ptsize{2.0pt};
      \draw[thick] (3.441458618106277,4.914912265733951) arc (55:125:6cm);
      \coordinate[label= above:\(A\)] (A) at (120:6cm);
      \coordinate[label=above:\(B\)] (B) at (60:6cm);
      \coordinate[label=above:\(C_1\)] (C) at (90:6cm);
      \draw (A) -- (B) -- (C) -- (A);
      \fill[fill=white,draw=black,thick] (C) circle (2pt);
      \fill[fill=white,draw=black,thick] (A) circle (2pt);
      \fill[fill=white,draw=black,thick] (B) circle (2pt);
      \node at (107.5:6.5cm) {$\gamma/2$};
      \node at (73.5:6.5cm) {$\gamma/2$};
    \end{tikzpicture}
    \caption{\label{AZ:fig:edge}%
      К доказательству Утв.~\ref{AZ:prep:1}
    }
  \end{figure}
  
  К хордам из
  второго и третьего множеств можно снова применить разбиение на три
  подмножества, рассмотрев середины дуг $AC_1$ и $C_1B$, и т.д. В результате
  получится покрытие счетным множеством <<вырожденных трубок>>: одна на дуге
  $\gamma/2$, две на дугах $\gamma/4$, четыре на дугах $\gamma/8$ и
  т.д. Суммарная мера, в силу счетной аддитивности и утверждения~\ref{AZ:prep:1}, равна
  $$
  \frac{\gamma^2}{2^2}+2\cdot \frac{\gamma^2}{4^2}+
  4\cdot \frac{\gamma^2}{8^2}+ 8\cdot \frac{\gamma^2}{16^2}+\ldots = \frac{\gamma^2}{2}\,.\qedhere
  $$
\end{proof}

\begin{prep}
  Мера Хаусдорфа $\bar{\mathcal H}^2(X)$ множества $X$ всех хорд равна
  $2\pi^2R^2$.
\end{prep}
\begin{proof}
  Разобъем огружность на  $n$ равных дуг. Тогда все множество хорд разбивается
  на $n(n-1)/2$ трубок, опирающихая на разные дуги, и $n$ множеств хорд с
  вершинами на одной дуге. Значит, из аддитивности меры и утверждений
  \ref{AZ:prep:tube}, \ref{AZ:prep:1} имеем
  $$
  n\cdot\frac{(2\pi R/n)^2}{2}+
  \frac{n(n-1)}{2}\cdot (2\pi R/n)^2=2\pi^2 R^2. \qedhere
  $$
\end{proof}

\bigskip

\textbf{3. Вероятностное пространство на множестве хорд окружности.}  %
Можно определить вероятностное пространство $(\Omega,\mathfrak F, \Pr)$, положив
$\Omega=X$ и выбрав в качестве сигма-алгебры $\mathfrak F$ семейство измеримых
относительно меры Хаусдорфа $\bar{\mathcal H}^2$ подмножеств множества $X$ и
положив
$$
\Pr(A)=\bar{\mathcal H}^2(A) / \bar{\mathcal H}^2(X).
$$

Будем называть хорду \emph{бертрановой}, если ее длина превосходит длину
$\sqrt{3} R$ вписанного в тот же круг правильного треугольника.

\begin{prep}
  Множество всех бертрановых хорд данного круга измеримо и его мера Хаусдорфа
  равна $(2\pi^2 R^2)/3$.
\end{prep}

\begin{proof}
  Чтобы доказать измеримость, построим монотонно возрастающую последователность
  измеримых множеств (т.е. из сигма-алгебры $\mathfrak F$), пределом которой
  является множество бертрановых хорд.

  \begin{figure}[h]
    \centering
      \begin{tikzpicture}
        \def\ptsize{2.0pt};
        \coordinate[label=above:$O$] (O) at (0,0);
        \draw[thick] (0,0) circle (2.5cm);
        \fill (O) circle(2pt);
        \draw[black!30,thick] (0:2.5cm)--(120:2.5cm)--(240:2.5cm)--cycle;
        \draw[black!30,thick] (-22.5:2.5cm)--(97.5:2.5cm)--(218.5:2.5cm)--cycle;
        \draw (0:2.5cm)--(135:2.5cm);
        \draw (-22.5:2.5cm)--(157.5:2.5cm);
        \draw[black!75] (0:2.5cm)--(157.5:2.5cm);
        \draw[black!75] (-22.5:2.5cm)--(180:2.5cm);
        \draw (0:2.5cm)--(180:2.5cm);
        \draw (-22.5:2.5cm)--(202.5:2.5cm);   
        \foreach \i in {0,...,16}
        \fill[fill=white,draw=black,thick] (22.5*\i:2.5cm) circle (2pt);
        \node at (0:2.5cm) [anchor=west] {$A_{k-1}$};
        \node at (-22.5:2.5cm) [anchor=west] {$A_{k}$};
        \node at (-22.5:2.5cm) [anchor=west] {$A_{k}$};
        \node at (135:2.5cm) [anchor=east] {$A^{(l)}$};
        \node at (157:2.5cm) [anchor=east] {$A^{(l-1)}$};
        \node at (168:2.6cm) [anchor=east] {$\cdot$};
        \node at (164:2.6cm) [anchor=east] {$\cdot$};
        \node at (171:2.6cm) [anchor=east] {$\cdot$};   
        \node at (202.5:2.5cm) [anchor=east] {$A'$};
        \node at (180:2.5cm) [anchor=east] {$A''$};
      \end{tikzpicture}
      \caption{\label{AZ:fig:cover}%
        К построению множества, включающего множество бертрановых хорд. Серыми
        линиями отмечены вспомогательные треугольники с вершинами в концах дуги
        $A_{k-1}A_k$}
   \end{figure}
   Разделим окружность на $n$ равных дуг точками $A_0$, $A_1$, \ldots, $A_{n-1}$
  (см. рис.~\ref{AZ:fig:cover}). Рассмотрим произвольную дугу $A_{k-1}A_k$. В
  каждой из точек $A_{k-1}$, $A_k$ восстановим вписанный треугольник. На
  противоположной части окружности (по отношению к дуге $A_{k-1}A_k$) эти
  треугольники ограниивают некоторое число $l$ точек $A'$, $A''$, \ldots,
  $A^{(l)}$. Конкретное значение $l$ мы установим позднее. Обозначим через
$$
    U_k=\bar{\mathcal T}(A_{k-1}A_k,A'A'')\cup
    \bar{\mathcal T}(A_{k-1}A_k,A''A''') \cup \ldots
    \cup \bar{\mathcal T}(A_{k-1}A_k,A^{(l-1)}A^{(l)})
$$
  объединение трубок, опирающихся на дугу $A_{k-1}A_k$ с одной стороны и на дуги
  $A'A''$, $A''A'''$, \ldots, $A^{(l-1)}A^{(l)}$ с другой стороны. Далее, повторим
  эту процедуру для всех дуг $A_{k-1}A_k$, $k=1$, $2$, \ldots, $n$ и рассмотрим множество
  $$
  V_n=\bigcup_{k=1}^n U_k.
  $$
  Подпоследовательность множеств $V_{2^n}$, $n\geqslant n_0$ можно сделать монотонной,
  если на кажом шаге делить пополам каждую из дуг. Тогда
  $V_{2^{n-1}}\subset V_{2^{n}}$ и  существует
  предел
  $$
  V=\lim_{n\to\infty} V_{2^n}=\bigcup_{n=n_0}^\infty V_{2^n},
  $$
  причем этот предел совпадает с множеством бертрановых дуг. Следовательно,
  $$
  \bar{\mathcal H}^{(2)}(V)=\lim_{n\to\infty} \bar{\mathcal H}^{(2)}(V_{2^n})
  $$
  
  Чтобы вычислить указанный предел, надо определить число $l$ (см. выше). Длина
  малой дуги равна $(2\pi R/n)$. Мы должны отбросить дуги, идущие в
  положительном направлении (против часовой стрелки) от точки $A_{k-1}$ до точки
  $A^{(l)}$, их всего $1+[(2\pi R/3)/(2\pi
  R/n)]=\nobreak\discretionary{}{\hbox{$\mathsurround=0pt =$}}{}1+[n/3]$ штук, и отбросить дуги в
  отрицательном направлении от $A_k$ до $A'$ в том же количестве. Остается
  $n-2([n/3]+1)$ дуг. Следовательно, множество $V_n$ есть объединение
  $n\cdot(n-2-2[n/3])/2$ различных трубок, каждая меры $(2\pi R/n)^2$. Значит,
  $$
  \bar{\mathcal H}^{(2)}(V_n)=\Bigl( \frac{2\pi R}{n}\Bigr)^2 \cdot \frac{n}{2} \cdot
  \Bigl(n-2-\Bigl[\frac{n}{3}\Bigr]\Bigl)\to \frac{2\pi^2 R^2}{3} \qquad
  \text{при $n\to\infty$.}
  $$
  Утверждение доказано.
\end{proof}

Как следствие, заключаем, что вероятность того, что случайная хорда бертранова,
равна $1/3$.

Интересно заметить, что мера на множестве хорд соответствует равномерному
распределению в треугольнике $\{(x,y)\colon 0\leqslant x<y<R\}$, где $x$,
$y$~--- длины дуг (в порядке возрастания) отсчитываемых от заданной точки. При
этом, диаметр треугольника (т.е. расстояние между наиболее удаленными друг от
друга точками) равен $R\sqrt{2}$, а не диаметру $2R$ множества хорд.

\end{document}